\documentclass[12pt,leqno]{article}
\usepackage{latexsym}

\newtheorem{defn}{Definition}[section]
\newtheorem{thm}[defn]{Theorem}
\newtheorem{cor}[defn]{Corollary}
\newtheorem{prop}[defn]{Proposition}

\def\A{\mbox{${\cal A}$}}
\def\F{\mbox{${\cal F}$}}
\def\U{\mbox{${\cal U}$}}
\def\O{\mbox{${\cal O}$}}
\newcommand{\rr}{\mbox{$I \hspace{-2.6mm} R$}}
\newcommand{\cc}{\mbox{$I \hspace{-2.6mm} C$}}
\newcommand{\zz}{\mbox{$I \hspace{-2.6mm} Z$}}
\newcommand{\nn}{\mbox{$I \hspace{-2.6mm} N$}}
\input{arrow.tex}

\begin{document}
\title{\bf On discontinuity of derivations, inducing non-unique complete metric topologies}

\author{S. R. PATEL}

\maketitle
\vspace{1.5cm}
\begin{tabbing}
\vspace{1.5cm}
\> Secondary 13F25, 46H40\\
\end{tabbing}

\noindent {\bf Abstract.} We give a simple method for constructing
commutative Fr\'{e}chet algebras which admit two inequivalent
Fr\'{e}chet algebra topologies. The result is applied to show that
the action of any non-algebraic analytic function may fail to be
uniquely defined among other useful applications. We give an
affirmative answer to a question of Loy from 1974. We also obtain
the uniqueness of the Fr\'{e}chet algebra topology of certain
Fr\'{e}chet algebras with finite dimensional radicals.
\newpage
\section{Introduction} The purpose of this paper is to give a simple
method for constructing commutative Fr\'{e}chet algebras which
have non-unique Fr\'{e}chet algebra topologies. In the Fr\'{e}chet
case, the only example of such an algebra so far known was that of
Read (see ~\cite{20}) an example originally constructed to show
the failure of the Singer-Wermer conjecture (the commutative
case), which holds in the Banach algebra case (see ~\cite{22}).
This algebra is $\F_\infty \,=\,\cc[[X_0,\, \dots, X_n, \dots]]$,
the algebra of all formal power series in infinitely many
commuting indeterminates $X_0,\, \dots, X_n, \dots$, and has a
Fr\'{e}chet algebra topology $\tau_R$ with respect to which the
natural derivation $\partial/\partial X_0$ is discontinuous with
the image the whole algebra. The other Fr\'{e}chet algebra
topology on $\F_\infty$ is the topology $\tau_c$ of coordinatewise
convergence with respect to which the derivation
$\partial/\partial X_0$ is continuous. Thus Read shows that the
situation on Fr\'{e}chet algebras is markedly different from that
on Banach algebras, and that a structure theory for Fr\'{e}chet
algebras behaves in a very distinctive manner from Banach
algebras. For example, the prestigious Michael problem is still
remained unsolved since 1952 (see ~\cite{11} for the recent
progress); Banach algebras are easily seen to be functionally
continuous.

Further we recall that in automatic continuity theory, we would
normally like to examine when and how the algebraic structure of
the algebra $A$ determines the topological structure of $A$, in
particular, the continuity aspect (and more particularly, the
uniqueness of the topology of $A$; see ~\cite{3, 5, 16, 17, 18}
for more details). So it is natural to expect that the
non-uniqueness of the topology would reflect some properties of
the algebraic structure of $A$. In ~\cite{13}, Feldman constructed
a Banach algebra with two inequivalent complete norms, in order to
show the failure of the Wedderburn Theorem. This algebra is
$\ell_2 \,\oplus \, \cc$ with one norm in which $\ell_2$ is dense.

It is easy to give uncountably many inequivalent Fr안chet space
topologies to a very familiar Fr안chet spaces, namely, spaces of
holomorphic functions ~\cite{23}. However the same space is also a
Fr안chet algebra of power series, and so, it admits a unique
Fr안chet algebra topology ~\cite{18}. Thus it is expected that the
case of having inequivalent Fr안chet algebra topologies by some
Fr안chet algebra should be rare, and the former situation should
be common.

We note that in ~\cite{15}, Loy gave a method for constructing
commutative Banach algebras with non-unique complete norm
topology, and we exploit the idea herein. We shall feel free to
use the terminology and conventions established there, and so our
argument here is kept short for the Fr\'{e}chet algebra case. We
use the discontinuity of derivations to give a Fr\'{e}chet algebra
$A$ other, inequivalent, Fr\'{e}chet algebra topologies. Although
the Fr\'{e}chet algebra topology $\tau_R\, + \,\tau_R$ of
$\F_\infty\;\oplus\;\F_\infty$ is not obtainable by our approach
(see below), other inequivalent Fr\'{e}chet algebra topologies on
$\F_\infty \oplus \F_\infty$ may be constructed. In view of the
use of the discontinuity of derivations method, it is interesting
to recall an early exposition of Banach algebras of powers series
and of their automorphisms and derivations, given by Grabiner in
~\cite{14}.

This simple method for constructing commutative Fr\'{e}chet
algebras with non-unique Fr\'{e}chet algebra topology surprisingly
also enables the demonstration of the uniqueness of the
Fr\'{e}chet algebra topology of certain Fr\'{e}chet algebras with
finite-dimensional radicals in Section 3. In the Fr\'{e}chet case,
though there are a few results on the uniqueness of the
Fr\'{e}chet algebra topology \linebreak (see ~\cite{3, 5, 11, 16,
18}), a class of Fr\'{e}chet algebras with finite-dimensional
radicals for this result appears to be treated for the first time
in this paper. Of course, in the Banach case, we have results on
the uniqueness of the complete norm topology of Banach algebras
with finite-dimensional radicals (see ~\cite{8, 19}); however, the
hypotheses and results are very different.

In Section 4, we obtain an affirmative answer to a question of Loy
from 1974 (see ~\cite{15}); we discuss the solution in the
Fr\'{e}chet case. In ~\cite{12}, Domar constructed a Banach
algebra with a quasinilpotent non-nilpotent radical; however, his
approach was different and used an entire function argument. In
particular, the question whether the construction of a Banach
algebra with a quasinilpotent non-nilpotent radical using the
method of the higher point derivations of infinite order is
possible, is a very good question. At present, we do not know the
answer. In the end, we discuss a particular example of a
semisimple Fr\'{e}chet algebra. Surprisingly, we note that the
converse of Theorem ~\ref{Therem 3_Loy} below does not hold.
\section{Non-uniqueness of the Fr\'{e}chet algebra topology} Let $A$ be a commutative metrizable LMC-algebra
and $M$ a commutative Fr\'{e}chet $A$-module, or simply
$A$-module, so that $M$ is a Fr\'{e}chet space which is a
commutative $A$-module such that the map $(x, m)\, \mapsto
\,x\cdot m$ is continuous from $A\,\times\,M$ to $M$. For such $A$
and $M$, let $H^1(A, M)$ denote the first algebraic cohomology
group, $H_C^1(A, M)$ the first continuous cohomology group, where
the cochains are required to be bounded. Thus with the usual
conventions $H^1(A, M)$ is the space of derivations of $A$ into
$M$, that is, linear mappings $D\,:\,A\,\rightarrow \,M$
satisfying $D(xy)\,=\,x\,\cdot\,D(y) \,+\,y\,\cdot\,D(x)$;
$H_C^1(A, M)$ the space of continuous derivations of $A$ into $M$.
As a simple example take $M\,=\,\cc$ with module action $x\,\cdot
\,\lambda\,=\lambda\,\phi(x)$ for some multiplicative linear
functional $\phi$ on $A$. If $\phi\,\neq\,0$, then $H^1(A, \cc)$
is the point derivation space at $\phi$.

Denoting the seminorms in both $A$ and $M$ by $(p_k)_{k\geq 1}$,
the constant $$l_k\,=\,\sup\{1,
\,p_k(xm)\,:\,p_k(x)\,=\,p_k(m)\,=\,1\}\;(k\,\in\,\nn),$$ is
finite, and if $\overline{A}$ is the completion of $A$ under
$(p_k)$, $M$ is clearly an $\overline{A}$-module. Let $\A$ denote
the vector space direct sum $\overline{A}\,\oplus\,M$ with product
$$(x, m) (y, n) \,=\,(xy, x\cdot n + y\cdot m)$$ and seminorms
$$q_k((x, m))\,=\,l_k(p_k(x) + p_k(m)).$$ For each $D\,\in\,H^1(A,
M)$, the functional $$q_{k, D}\,:\,(x, m)\,\rightarrow\,l_k(p_k(x)
+ p_k(D(x)-m))$$ is defined on the subalgebra $A\,\oplus\,M$ of
$\A$ and is easily seen to be a submultiplicative seminorm
thereon. If $D\,=\,0$ the completion of $A\,\oplus\,M$ under
$(q_{k,D})$ is of course just $\A$. For general $D$ we note that
the map
$$\theta_D\,:\,(x, m)\,\rightarrow\,(x, D(x) - m)$$ is an isometric
isomorphism of $A\,\oplus\,M$ under $(q_{k,D})$ into $\A$, and so
extends uniquely to a map of the completion $\A_{D}$ of
$A\,\oplus\,M$ into $\A$. In particular, if
$\iota\,:\,x\,\rightarrow\,(x, 0)$ is the natural embedding of $A$
into $A\,\oplus\,M$ then $q_{k, D}^{'}\,:\,x\,\rightarrow\,q_{k,
D}(\iota(x))$ is a seminorm on $A$ and $\theta_D\, \circ \,\iota$
extends to an isometric isomorphism of $\overline{A_D}$, the
completion of $A$ under $q_{k, D}^{'}$, with
$\overline{\textrm{Gr}(D)}$, the closure (in $\A$) of the graph of
$D$.

Now if $D$ is continuous, then $(q_{k, D})$ is equivalent to
$(q_k)$ on $A\,\oplus\,M$ and $(q_{k, D}^{'})$ is equivalent to
$(p_k)$ on $A$. Thus $\A_{D}\,=\,\A$ and
$\overline{A_D}\,=\,\overline{A}$. In the discontinuous case
$q_{k, D}^{'}$ is a discontinuous seminorm on $A$ and $\iota$ is a
discontinuous isomorphism. This latter result has been used in the
Fr\'{e}chet case as follows: let $A$ be the algebra of polynomials
on a fixed open neighbourhood $U$ of the closed unit disc $\Delta$
with the compact-open topology, $M\,=\,\cc$ with module action
$p\cdot \lambda\,=\,\lambda\,p(1)$ for
$(p,\,\lambda)\,\in\,A\,\oplus\,M$ and
$D\,:\,p\,\rightarrow\,p'(1)$. Here $\overline{A_D}$ is a singly
generated Fr\'{e}chet algebra with spectrum $U$ and one
dimensional radical (for the Banach case, see ~\cite{10, 15}).

Suppose now that $D$ is discontinuous. If $M$ is finite
dimensional, then it easily follows that
$\overline{\textrm{Gr}(D)}\,\bigcap\,\{0\}\,\oplus\,M\,\neq\,\{0\}$,
and if $A$ is complete this holds for arbitrary $M$ by the closed
graph theorem for Fr\'{e}chet spaces. Thus if $A$ is a Fr\'{e}chet
algebra with $H^1(A, M)\,\neq\,H_C^1(A, M)$ for some $M$, then $A$
has a completion with a non-trivial nil ideal. For example,
$A\,=\,\F_\infty$ is a Fr\'{e}chet algebra under the Fr\'{e}chet
algebra topology $\tau_R$ imposed by Read, then $H^1(A,
M)\,\neq\,H_C^1(A, M)$ for $M\,=\,\F_\infty$, since the derivation
$\partial/\partial X_0$ is discontinuous with respect to this
topology. Thus $\F_\infty$ has a completion with a non-trivial nil
ideal. We remark that if $A$ is a semisimple Fr\'{e}chet algebra,
then $H^1(A, M)\,=\,H_C^1(A, M)$ ~\cite{4}. However we do not know
an example of a (semisimple) non-Banach Fr\'{e}chet algebra such
that $H^1(A, M)\,=\,0$ for any $A$-module $M$; the Singer-Wermer
conjecture holds for commutative semisimple Banach algebras.
Banach algebras with $H^1(A, M)\,=\,H_C^1(A, M)$ are discussed in
~\cite{2}. Finally we add one further hypothesis to these
considerations to obtain the following result.
\begin{thm} \label{Theorem 1_Loy} Let $A$ be a commutative Fr\'{e}chet
algebra, $D$ a non-zero derivation of $A$ into a commutative
Fr\'{e}chet $A$-module $M$. If $D$ vanishes on a dense \linebreak
subset of $A$ then the algebra $\overline{A_D}$ admits two
inequivalent Fr\'{e}chet algebra topologies.
\end{thm}
{\it Proof.} The proof is the same as that of ~\cite[Theorem
1]{15}. $\hfill \Box$

As a corollary, we have the following special case
\begin{cor} \label{Corollary 1_Read} Let $A$ be the Fr\'{e}chet algebra $\F_\infty$
under the Fr\'{e}chet algebra topology $\tau_R$ and let $D$ be the
natural derivation $\partial/\partial X_0$. Then
$\overline{(\F_\infty)}_D$ admits another Fr\'{e}chet algebra
topology $\tau_D$, generated by $(q_{k, D})$, different from
$\tau_R + \tau_R$, generated by $(q_k)$.
\end{cor}
{\it Proof.} By ~\cite[Theorem 2.5]{20}, $\partial/\partial X_0$
is a discontinuous derivation on $(\F_\infty, \tau_R)$ which
vanishes on a dense subset of $(\F_\infty, \tau_R)$ since
$X_n\,\rightarrow \,X_0$ in $(\F_\infty, \tau_R)$; $X_0$ lies in
the closure of the coefficient algebra $\A_0\,=\,\cc[[X_1,\,\dots,
\,X_n,\,\dots]]$ of $\F_\infty\,=\,\cc[[X_0,\,\dots,
\,X_n,\,\dots]]$. Thus, by Theorem ~\ref{Theorem 1_Loy},
$\F_\infty\,\oplus\,\F_\infty$ is also a Fr\'{e}chet algebra under
$\tau_D$, generated by $(q_{k, D})$, different from $\tau_R +
\tau_R$, generated by $(q_k)$. $\hfill \Box$

In fact, we have a more general result than Theorem ~\ref{Theorem
1_Loy} as follows. Let $A$ be a commutative Fr\'{e}chet algebra,
and let $M$ be a Fr\'{e}chet $A$-module. Set $\U
\,=\,A\,\oplus\,M$, where $(a, x) (b, y)\,=\,(ab, a\cdot y +
b\cdot x)$ for $a, b\, \in\, A$ and $x, y \,\in \,M$. Then $\U$ is
a commutative algebra with $\textrm{Rad}\; \U\,=\,\textrm{Rad}\;
A\,\oplus\,M$. Let $D\,:\,A\,\rightarrow\,M$ be a derivation, and
set
$$q_k((a, x))=p_k(a)+p_k(x), \;q_{k, D}((a, x))=p_k(a)+p_k(D(a) -
x) \;(a \in A,\,x \in M).$$
\begin{thm} \label{Thm. 5.1.17_Frechet} The algebra $\U$ is a Fr\'{e}chet
algebra with respect to both $(q_k)$ and $(q_{k, D})$. The two topologies
are equivalent if and only if $D$ is continuous.
\end{thm}
{\it Proof.} Certainly $(\U, (q_k))$ is a Fr\'{e}chet algebra and
$ q_{k, D}$ is a seminorm on $\U$ for each $k\,\in\,\nn$. For $(a,
x),\, (b, y)\,\in\,\U$, we have $q_{k, D}((a, x)(b,
y))\,=\,p_k(ab)+p_k(a \cdot (D(b) - y) + b \cdot (D(a) - x)) \leq
(p_k(a) + p_k(D(a) - x)) (p_k(b) + p_k(D(b) - y)) = q_{k, D}((a,
x)) q_{k, D}((b, y))$, and so $ q_{k, D}$ is a submultiplicative
seminorm on $\U$ for each $k\,\in\,\nn$. We now show that $(\U,\,
(q_{k, D}))$ is a Fr\'{e}chet algebra. Let $((a_n, x_n))$ be a
Cauchy sequence in $(\U,\, (q_{k, D}))$. Then $(a_n)$ and $(D(a_n)
- x_n)$ are Cauchy sequences in $(A, \,(p_k))$ and $(M,\,(p_k))$,
respectively. Since $A$ and $M$ are Fr\'{e}chet spaces, there
exists $a\,\in\,A$ and $x\,\in\,M$ such that $a_n\,\rightarrow\,a$
and $D(a_n) - x_n\,\rightarrow\,x$. Then $(a_n,
x_n)\,\rightarrow\,(a, D(a) - x))$ in $(\U,\, (q_{k, D}))$ and so
$(\U,\, (q_{k, D}))$ is a Fr\'{e}chet algebra.

Suppose that $D$ is continuous. Then, for each $m\,\in\,\nn$,
there exists $n(m)\,\in\,\nn$ and a constant $c_m\,>\,0$ such that
$$q_{m, D}((a, x))\leq p_{m}(a)+c_mp_{n(m)}(a)+p_m(x)\leq(1+c_m)q_{n(m)}((a,
x))\;((a, x)\in \U),$$ and so the two topologies are equivalent,
by the open mapping theorem for Fr\'{e}chet spaces.

Conversely, suppose that the two topologies are equivalent on the
algebra $\U$. Then, for each $m\,\in\,\nn$, there exists
$n(m)\,\in\,\nn$ and a constant $c_m\,>\,0$ such that $ q_{m,
D}((a, x))\,\leq\,c_m\,q_{n(m)}((a, x))\;\;((a, x)\,\in\,\U)$.
Hence
$$p_m(D(a))\,\leq\, q_{m, D}((a, 0))\,\leq\,c_m\,q_{n(m)}((a,
0))\,=\,c_m\,p_{n(m)}(a)\;\;(a\,\in\,A),$$ and so $D$ is
continuous. $\hfill \Box$

As corollaries, we have the following important results.
\begin{cor} \label{Thm. 5.1.18_Frechet} There is a commutative algebra
with a one-dimensional radical which is a Fr\'{e}chet algebra with
respect to two inequivalent Fr\'{e}chet algebra topologies.
\end{cor}
{\it Proof.} Let $A$ be a Fr\'{e}chet function algebra with a
discontinuous point derivation $D$ at a continuous character
$\phi$. Then $\cc$ is a Fr\'{e}chet $A$-module with respect to the
operation $(f,
z)\,\mapsto\,\phi(f)z,\;\;A\,\times\,\cc\,\rightarrow \,\cc$, and
so we are in a situation where Theorem ~\ref{Thm. 5.1.17_Frechet}
applies: $\U\,=\,A\,\oplus\,\cc$ is a Fr\'{e}chet algebra with
respect to the product $$(f, z) (g, w)\,=\,(fg,
\phi(f)w\,+\,\phi(g)z)$$ and each of the topologies generated by
$(q_k)$ and $(q_{k, D})$, respectively, where $q_k((f,
z))\,=\,p_k(f)\,+\,\mid z\mid$ and $q_{k, D}((f,
z))\,=\,p_k(f)\,+\,\mid D(f) - z\mid$.

As above, $\textrm{Rad}\;\U\,=\,\{0\}\,\oplus\,\cc$, and so
$\textrm{Rad}\;\U$ is one-dimensional. $\hfill \Box$

\begin{cor} \label{Corollary 2_Read} Let $A$ be the Fr\'{e}chet
algebra $\F_\infty$ under the Fr\'{e}chet algebra topology $\tau_c$
and let $D$ be the natural derivation $\partial/\partial X_0$. Then
the two topologies $\tau_D$, generated by $(q_{k, D})$, and $\tau_c + \tau_c$,
generated by $(q_k)$, are equivalent on the algebra $\F_\infty\,\oplus\,\F_\infty$
if and only if $D$ is continuous.
\end{cor}
{\it Proof.} By ~\cite{20}, $\partial/\partial X_0$ is a
continuous derivation on $(\F_\infty, \tau_c)$. Thus, by Theorem
~\ref{Thm. 5.1.17_Frechet}, the two topologies $\tau_D$, generated
by $(q_{k, D})$, and $\tau_c + \tau_c$, generated by $(q_k)$, are
equivalent on the algebra $\F_\infty\,\oplus\,\F_\infty$. Converse
is trivial. $\hfill \Box$

As an application of this method, let $A$ be the algebra
$\textrm{Hol}(U)$ of analytic functions on the open unit disc $U$,
with the compact-open topology. Let $\O$ be the algebra of
functions analytic in a neighbourhood of $U$, and $\psi : A
\rightarrow \O$ the \linebreak monomorphism $\psi x(\lambda) =
x(\lambda/2),\;x \in A, \;\mid \lambda \mid \,< 2$. Finally let
$M$ be a Banach space, $T$ an endomorphism of $M$ with norm at
most $1$. Then, following the argument given in ~\cite{15}, $M$ is
an $A$-module. In particular, we may consider the situation
developed in ~\cite{6} but in the Fr\'{e}chet case, where
$M\,=\,C[0, 1]$ and $T$ is the operator of indefinite integration.
Non-zero linear mappings $\beta : \O \rightarrow M$ are
constructed in ~\cite{6} to vanish on polynomials and satisfy
$$\beta(fg) \,=\,f(T)\beta(g)\,+\,g(T)\beta(f)\;\;(f,
g\,\in\,\O).$$ Letting $D\,=\,\beta\psi$ we have that
$D\,:\,A\,\rightarrow\,M$ is a derivation vanishing on
polynomials. Since $D\,\neq\,0$, it is necessarily discontinuous
(which answers a question of ~\cite{2} in the Fr\'{e}chet case).

We also note that $A$ and $D$ here satisfy the hypothesis of
Theorem ~\ref{Theorem 1_Loy} since polynomials are dense in $A$.
Thus $\overline{A_D}$ has two inequivalent Fr\'{e}chet algebra
topologies with infinite dimensional radical (by the construction
of $\beta$; see ~\cite{6}). Indeed algebras with finite
dimensional radical have unique functional calculus by the
Fr\'{e}chet analogue of Theorem 1 of ~\cite{6} so that the
argument above shows that any derivation of $A$ into a finite
dimensional $A$-module which vanishes on polynomials must be zero.
In fact, such a derivation is continuous, since $A$ is an algebra
of class $\mathcal{B}$ (see ~\cite{9}) with $A\,\oplus\,M$ a
strongly decomposable Fr\'{e}chet algebra with finite dimensional
radical $M$, and so $\theta_D\,\iota$ is continuous by Theorem
~\ref{Theorem 5.4_DM} below.

\noindent {\bf Remark.} The results established by Dales in
~\cite{6} have appropriate analogues in the Fr\'{e}chet case.

Using the present example we can show that the exponential
function in such an algebra is not independent of the Fr\'{e}chet
algebra topology (see ~\cite[p. 413]{15} for details).

We recall that the first discontinuous functional calculus map was
constructed by Allan (see ~\cite[Theorem 8]{1}). The algebras
satisfying this theorem are $L_{\textrm{loc}}^1(\rr^+)$,
$L^1(\rr^+, W)$ and $C^\infty(\rr^+)$. It follows from proof of
the theorem that the above condition is sufficient on an algebra
$A$ for the existence of a discontinuous functional calculus
homomorphism. However, the algebra
$\U\,=\,\textrm{Hol}(U)\,\oplus\,\cc$ shows that it is not a
necessary condition. In fact, it is shown in Theorem 8 of
~\cite{1} that there is an incomplete metrizable topology on
$\textrm{Hol}(U)$ which dominates the compact-open topology, and
the completion of $\textrm{Hol}(U)$ in such a topology has a
nilpotent radical, as discussed before. Such a result is
impossible for $\F \,=\,\cc[[X]]$ (since $X$ is locally nilpotent)
and $C(U)$ (no adequate theory of point derivations).
\section{Uniqueness of the Fr\'{e}chet algebra topology} In the above situations the ideal adjoined was
always nilpotent of index two; we now consider how to obtain more
general ideals by following ~\cite[\S 2]{15}.

Let $A$ be a commutative metrizable LMC-algebra, $M$ an $A$-module
which is also a commutative Fr\'{e}chet algebra, $D =
\{D_1,\,\dots,\, D_r\}$ a higher derivation of rank $r$ of $A$
into $M$. In analogy to $\A$ let $\A_r$ denote the vector space
$\overline{A}\,\oplus\,M^r$ with product as given in ~\cite[p.
414]{15}, and seminorms
$$q_k((x,\;\{m_i\}))\;=\;l_k(p_k(x)\;+\;\sum_{i=1}^{r}p_k(m_i)).$$ We also have the seminorms
$q_{k, D}$ on $A\;\oplus \;M^r$,
$$q_{k, D}((x,\;\{m_i\}))\;=\;l_k(p_k(x)\;+\;\sum_{i=1}^{r}p_k(D_i(x) - m_i)),$$ the isomorphism
$\theta_D\;:\;(x,\;\{m_i\})\;\rightarrow\;(x,\;\{D_i(x) - m_i\})$
and the completion $\overline{A_D}$ of $A$ under $q_{k,
D}^{'}\;:\;x\;\mapsto\;q_{k, D}((\iota(x))$. We consider a
specific case below.

Thus let $A$ be the algebra of polynomials on a fixed open
neighbourhood $U$ of the closed unit disc $\Delta$, with seminorms
$p_k$ generating the compact-open topology, $M\;=\;\cc$ with
module action as before and a higher point derivation $D$ of rank
$r$ as given in ~\cite[p. 415]{15}, so that $$q_{k,
D}{'}(p)\;=\;p_k(p)\;+\;\sum_{i=1}^{r}\frac{|p^{(i)}(1)|}{i!}.$$
Then, the arguments given in Section 2 of ~\cite{15} for the
Banach case can also be applied to the Fr\'{e}chet case to see
that $\overline{A_D}$ has a radical which is nilpotent of index
$r$.

We note that in this example $\A_r$, and hence $\overline{A_D}$,
is a strongly decomposable Fr\'{e}chet algebra of class
$\mathcal{B}$ with finite dimensional radical and so has a unique
Fr\'{e}chet algebra topology by Corollary ~\ref{Corollary 5.6_DM}
below. We remark that Section 5 of ~\cite{9} could be extended to
the Fr\'{e}chet case (up to Proposition 5.5). In particular, we
have the following
\begin{thm} \label{Theorem 5.4_DM} Let $B$ be an algebra of class
$\mathcal{B}$ and let $A$ be a strongly decomposable Fr\'{e}chet
algebra with finite dimensional radical $R$. Then any homomorphism
$\theta\;:\;B\;\rightarrow\;A$ is continuous. $\hfill \Box$
\end{thm}
\begin{prop} \label{Proposition 5.5_DM} Let $A$ be an algebra of class
$\mathcal{B}$ with finite dimensional radical $R$. Then any
decomposition of $A$ is a strong decomposition. $\hfill \Box$
\end{prop}

Then, as a corollary to Theorem \ref{Theorem 5.4_DM}, we have the
following result.
\begin{cor} \label{Corollary 5.6_DM} If $A$ is a decomposable Fr\'{e}chet
algebra of class $\mathcal{B}$ with finite dimensional radical
$R$, then $A$ has a unique Fr\'{e}chet algebra topology. $\hfill
\Box$
\end{cor}

In comparison to this we note the following result the proof of
which is omitted.
\begin{thm} \label{Therem 2_Loy} Let $A$ be a commutative Fr\'{e}chet
algebra, $D\;=\;\{D_1,\;\dots,\;D_r\}$ a higher point derivation
of rank $r$ of $A$ into $\cc$ such that $\{D_1,\;\dots,\;D_r\}$ is
a set of discontinuous functionals. Then $\overline{A_D}$ admits
two inequivalent Fr\'{e}chet algebra topologies and has nilpotent
elements of index $r$. $\hfill \Box$
\end{thm}
\section{Affirmative answer to Loy's question}
In ~\cite{15}, Loy raised the question of whether quasinilpotent
non-nilpotent radicals are obtainable in some analogous fashion.
We now answer this question for a more general case of Fr\'{e}chet
algebras as follows. We remark that the Jacobson radical in a
commutative Fr\'{e}chet algebra $A$ may also be defined as the set
of quasinilpotent elements, that is, $\textrm{Rad}
A\,=\,\{x\,\in\,A\,:\,r(x) \,=\,0\}$, where $r(x)
\,=\,\sup_{k\,\in\,\nn}r_k(x_k)$.

First, we consider the Banach algebra situation. So let $A$ be a
commutative normed algebra, $M$ an $A$-module which is also a
commutative Banach algebra, $D\,=\,\{D_1,\,\dots\}$ a higher
derivation of infinite order of $A$ into $M$, so that for $x,
y\,\in\,A$ and $s\,\in\,\nn$,
$$D_s(xy)\,=\,x \, D_s(y) + y \, D_s(x) + \sum_{i =1}^{s
- 1}D_i(x) D_{s -  i}(y).$$ In analogy to $\A_r$ let $\A_\infty$
denote the vector space $\overline{A}\,\oplus\,M^\infty$ with
product $$(x,\;\{m_i\})\;(y,
\;\{n_i\})\;=\;(xy,\;\{xn_i\;+\;ym_i\;+\;\sum_{j=1}^{i-1}m_jn_{i-j}\})$$
and metric
$$d((x,\;\{m_i\}), 0)\;=\;\|x\|\;+\;\sum_{i=1}^{\infty}\frac{2^{-i}\|m_i\|}{1+\|m_i\|}.$$
We also have the metric $d_D$ on $A\;\oplus \;M^\infty$,
$$d((x,\;\{m_i\}), 0)_D\;=\;\|x\|\;+\;\sum_{i=1}^{\infty}\frac{2^{-i}\|D_i(x) - m_i\|}{1+\|D_i(x) - m_i\|},$$
the isomorphism
$\theta_D\;:\;(x,\;\{m_i\})\;\rightarrow\;(x,\;\{D_i(x) - m_i\})$
and the completion $\overline{A_D}$ of $A$ under the metric $x
\mapsto d(\iota(x), 0)_D$. We consider a specific case below.

Thus let $A$ be the algebra of polynomials on the closed unit disc
$\Delta$, with the uniform norm $\|\cdot\|_\infty$, $M \,= \,\cc$
with module action as before and a higher point derivation $D$ of
infinite order as given in ~\cite[p. 415]{15}, so $M^\infty =
\cc_0[[X]]$, with powers in $M^\infty$ move to the right due to
the convolution product (and thus, they would eventually be zero
in $M^r, r \in \nn$); of course, continuity of multiplication is
apparent if one thinks of the usual coordinatewise \linebreak
convergence topology $\tau_c$ on $\cc[[X]]$, which is equivalent
to $d$. Following the arguments given on p. 415 of ~\cite{15}, we
wish to show that given a sequence $(\alpha_i)$ \linebreak of
complex numbers, there exists a sequence $(p_n)$ in $A$ with
$\|p_n\|_\infty \rightarrow 0$ and\linebreak $D_i(p_n) \rightarrow
\alpha_i$ for each $i$. By ~\cite[p. 415]{15}, for each fixed $k
\in \nn$, there is a sequence $(p_n^k)$ in $A$ such that
$\|p_n^k\|_\infty \to 0$ and $D_i(p_n^k) \to \alpha_i$ for $i =
1,\,\ldots,\,k$. For each $k \in \nn$, choose a polynomial $p_k =
p_n^k$ from the sequence $(p_n^k)$ by taking $n$ sufficiently
large such that $\|p_k\|_\infty\,<\,\frac{1}{k}$ and $|D_i(p_k) -
\alpha_i|\,<\,\frac{1}{k}$ for $i = 1,\,\ldots,\,k$. Then we have
a sequence $(p_k)$ such that $\|p_k\|_\infty\,<\,\frac{1}{k}$ and
for $i \in \nn$, for any $k > i$, $|D_i(p_k) - \alpha_i|
\,<\,\frac{1}{k}$. Thus $\|p_k\|_\infty \to 0$ and for every $i
\in \nn$, $D_i(p_k) \to \alpha_i$. The existence of the required
sequence\linebreak $(p_n)$ is now clear. Now let $\Theta$ denote
the extension by continuity of the isomorphism $\theta_D\iota$ to
$\overline{A_D}$. Then if $(x, \{\alpha_i\}) \in A \oplus
M^\infty$ we have $(x, \{\alpha_i\}) = (x, \{D_ix\}) + (0,
\{\alpha_i - D_ix\}) \in \Theta(\overline{A_D})$, and since
$\Theta $ is an isometry and $\overline{A_D}$ is complete, it
follows that $\Theta$ maps $\overline{A_D}$ onto $\A_\infty =
\overline{A} \oplus M^\infty$ which contains the \linebreak
radical $\{0\} \oplus M^\infty$ (note that, by an analogous
argument, preceding to Theorem 2.3, $\textrm{Rad}\A_\infty =
\textrm{Rad}\overline{A} \oplus M^\infty = \{0\} \oplus
M^\infty$). Thus $\overline{A_D}$ has a radical which has
quasinilpotent non-nilpotent elements. We do not know whether
$\overline{A_D}$ has a unique Fr\'{e}chet algebra topology.
However we have the following result the proof of which is omitted
(see ~\cite[2.2.46 (ii)]{7} for details on the Dales-McClure
Banach algebra, and to know about the two inequivalent Fr\'{e}chet
algebra topologies on $\overline{A_D}$, follow either Theorem
~\ref{Therem 2_Loy} or Theorem 2.3).
\begin{thm} \label{Theorem 2_Loy} Let $A$ be the Dales-McClure Banach
algebra, $D\,=\,(D_i)$ a totally discontinuous higher point
derivation of infinite order at a character $\phi$. Then
$\overline{A_D}$ admits two inequivalent Fr\'{e}chet algebra
topologies and has quasinilpotent non-nilpotent elements. $\hfill
\Box$
\end{thm}

Now we can consider the Fr\'{e}chet case by combining the
situations given in \S 3 and in the Banach case above. We consider
a specific case below.

Thus let $A$ be the algebra of polynomials on a fixed open
neighbourhood $U$ of the closed unit disc $\Delta$, with seminorms
$p_k$ generating the compact-open topology, $M\;=\;\cc$ with
module action as before and a higher point derivation $D$ of
infinite order as above, so that
$$q_{k, D}{'}(p)\;=\;p_k(p)\;+\;
\sum_{i=1}^{\infty}\frac{|p^{(i)}(1)|}{i!}.$$ Following the
arguments given in the Banach case above, we see that
$\overline{A_D}$ has a radical which has quasinilpotent
non-nilpotent elements. We do not know whether $\overline{A_D}$
has a unique Fr\'{e}chet algebra topology. However we have the
following result the proof of which is omitted. We note that the
Dales-McClure Fr\'{e}chet algebra can be constructed along the
lines of the Dales-McClure Banach algebra by replacing a weight
$\omega$ on $\zz^+$ by an increasing sequence $W\, =\, (\omega_k)$
of weights on $\zz^+$. Thus,
$\check{V}_WE\;=\;\bigcap_{k=1}^{\infty}\check{V}_{\omega_k}E$, an
analogue of the Beurling-Fr\'{e}chet algebra (see ~\cite[Example
1.2]{3}).
\begin{thm} \label{Theorem 2_Loy} Let $A$ be the Dales-McClure Fr\'{e}chet
algebra, $D\,=\,(D_i)$ a totally discontinuous higher point
derivation of infinite order at a continuous character $\phi$.
Then $\overline{A_D}$ admits two inequivalent Fr\'{e}chet algebra
topologies and has quasinilpotent non-nilpotent elements. $\hfill
\Box$
\end{thm}

We now obtain some special results when $A$ is the algebra of
polynomials on a fixed open neighbourhood $U$ of the closed unit
disc $\Delta$. Let $(M_k)$ be a sequence of positive reals such
that
$$\frac{M_k}{k!}\;\geq\;\frac{M_i}{i!}\frac{M_{k-i}}{(k-i)!}$$ for $1\;\leq\;i\;<\;k$
so that, in particular, the sequence
$k\;\mapsto\;(\frac{k!}{M_k})^{\frac{1}{k}}$ is monotonic
decreasing, with limit $\gamma\;\geq\;0$. Consider the seminorms
$q_k$ on $A$ given by $$q_k(p)
\;=\;p_k(p)\;+\;\sum_{i=1}^{\infty}\frac{|p^{(i)}(1)|}{M_i}$$ and
let $\overline{A_\infty}$ be the completion of $A$ under $(q_k)$.
Then we have the following result the proof of which is omitted
(see ~\cite[Theorem 3]{15} for details in the Banach case).
\begin{thm} \label{Therem 3_Loy} $\overline{A_\infty}$ is
semisimple if $\gamma \; >\;0$. For each $k\;\in\;\nn$, $(A, q_k)$
is natural if and only if $\gamma\;=\;0$. $\hfill \Box$
\end{thm}

In fact, if $\overline{A_\infty}\,=
\,\lim\limits_{\longleftarrow}(\overline{A_\infty})_{k}$ is the
Arens-Michael representation, then each
$(\overline{A_\infty})_{k}$ satisfies Theorem 3 of ~\cite{15} (and
hence, the use of $(M_k)$ is implicit and necessary for the
Fr\'{e}chet case). If $\gamma = 0$, then each $(A, q_k)$ is
natural and vice-versa. In this case, the inverse limit algebra
$(A, (q_k))$ is also natural. In the case $\gamma > 0$, each $(A,
q_k)$ is not natural and the completion
$(\overline{A_\infty})_{k}$ is semisimple, and so,
$\overline{A_\infty} \,=\,
\lim\limits_{\longleftarrow}(\overline{A_\infty})_{k}$ is also a
semisimple Fr\'{e}chet algebra. However we note that for the
inverse limit algebra $(A, (q_k))$ to be natural, then it is not
necessary that each $(A, q_k)$ is natural. We exhibit this
interesting case as follows. For each fixed $k\in \nn$, let
$\{M_i^k\}$ be a sequence of positive reals with properties as
given above. Consider the seminorms $q_k^{'}$ on $A$ given by
\linebreak
$$q_k^{'}(p)
\;=\;p_k(p)\;+\;\sum_{i=1}^{\infty}\frac{|p^{(i)}(1)|}{M_i^k}$$
and let $\overline{A_\infty}$ be the completion of $A$ under
$(q_k^{'})$. We select the sequence $\{M_i^k\}$ of sequences of
positive reals such that the limit $\gamma_k\;>\;0$ for each fixed
$k$, but $\gamma_k\;\rightarrow\;0$ as $k\;\rightarrow\;\infty$,
that is, $\gamma\;=\;\lim\gamma_k\;=\;0$. Thus, by ~\cite[Theorem
3]{15}, each $(A, q_k^{'})$ is not natural and the completion
$(\overline{A_\infty})_{k}$ is semisimple. Then
$\overline{A_\infty}\,=
\,\lim\limits_{\longleftarrow}(\overline{A_\infty})_{k}$ is a
semisimple Fr\'{e}chet algebra. This shows that the converse of
Theorem ~\ref{Therem 3_Loy} does not hold. Also the inverse limit
algebra $(A, (q_k^{'}))$ is natural. To exhibit the above
situation, take $M_i^k\;=\;k^ii!$ so that
$\gamma_k\;=\;\frac{1}{k}$ with $\lim\gamma_k\;=\;0$. On the other
hand, Rolewicz in ~\cite{21} constructed an example of a
semisimple Fr\'{e}chet algebra which is not an inverse limit of
semisimple Banach algebras. So it is possible to have
$\overline{A_\infty}\,=
\,\lim\limits_{\longleftarrow}(\overline{A_\infty})_{k}$ a
semisimple Fr\'{e}chet algebra such that each
$(\overline{A_\infty})_{k}$ is not a semisimple Banach algebra,
and so, by Theorem 3 of ~\cite{15}, $\gamma_k\;=\;0$ for each $k$
(with $\gamma\;=\;\lim\gamma_k\;=\;0$) which implies that each
$(A,\;q_k^{'})$ is natural and the inverse limit algebra
$(A,\;(q_k^{'}))$ is also natural.

\noindent {\bf Acknowledgment.} The author thanks Professor
Richard J Loy of the Australian National University (Canberra) for
fruitful comments, especially for his valuable comments pertaining
to my attempt to his question, above.

{\bf Address:} Ahmedabad, Gujarat, INDIA. E-mails:
srpatel.math@gmail.com, coolpatel1@yahoo.com
\end{document}